\documentclass{amsart}

\usepackage{hyperref}
\hypersetup{
    colorlinks=true,
    linkcolor=red,
    urlcolor=pink,
    citecolor=blue
}
\usepackage{amssymb}
\usepackage{amsmath}
\usepackage{amsthm}
\usepackage{enumerate}
\usepackage{graphicx}
\usepackage{color}
\usepackage{subcaption}
\usepackage{float}

\usepackage{tikz}
\usetikzlibrary{arrows.meta}

\usepackage[english]{babel}
\usepackage[utf8]{inputenc}

\usepackage{multicol}
\usepackage{xcolor}
\usepackage{array}

\theoremstyle{plain}

\makeatletter
\newtheorem*{rep@theorem}{\rep@title}
\newcommand{\newreptheorem}[2]{%
\newenvironment{rep#1}[1]{%
 \def\rep@title{#2 \ref{##1}}%
 \begin{rep@theorem}}%
 {\end{rep@theorem}}}
\makeatother

\newreptheorem{theorem}{Theorem}

\theoremstyle{definition}

\theoremstyle{remark}

\numberwithin{equation}{section}

\newtheorem{thm}{Theorem}[section]

\newtheorem{cor}[thm]{Corollary}
\newtheorem{con}[thm]{Conjecture}
\newtheorem{ques}[thm]{Question}

\newcommand{\mC}{\mathbb{C}}
\newcommand{\mR}{\mathbb{R}}

\newcommand{\mD}{\mathbb{D}}

\newcommand{\hC}{\widehat{\mathbb{C}}}
\newcommand{\greencheck}{\textcolor{green!70!black}{\ensuremath{\checkmark}}}
\newcommand{\redcross}{\textcolor{red}{\ensuremath{\boldsymbol{\times}}}}
\DeclareMathOperator{\dist}{dist}

\DeclareMathOperator{\Mod}{Mod}
\DeclareMathOperator{\caplog}{cap}
\DeclareMathOperator{\diam}{diam}

\definecolor{domainpink}{RGB}{250,190,190}
\tikzset{
  hole/.style={draw=black, line width=.43pt, fill=white},
  transformation/.style={-{Stealth[length=2mm]}, line width=.55pt}
}

\begin{document}

\title{Conformal rigidity: progress and challenges
}

\date{September 30, 2026}

\author[D. Ntalampekos]{Dimitrios Ntalampekos}
\address{Department of Mathematics, Aristotle University of Thessaloniki, Thessaloniki, Greece.}
\email{dntalam@math.auth.gr}

\author[M. Younsi]{Malik Younsi}
\address{Department of Mathematics, University of Hawaii Manoa, Honolulu, USA.}
\email{malik.younsi@gmail.com}

\thanks{D.N. was supported by the European Union (ERC, GRComPaS, Grant Agreement no.\ 101214615). M.Y. was supported by NSF Grant DMS-2350530.}

\begin{abstract}
A domain in the Riemann sphere is called a \textit{circle domain} if every connected component of its boundary is either a round circle or a single point. A circle domain is \textit{conformally rigid} if every conformal map onto another circle domain is the restriction of a M\"{o}bius transformation. Conformal rigidity is closely related to Koebe's conjecture, which asserts that every domain in the Riemann sphere is conformally equivalent to a circle domain. In this article, we survey recent progress on conformal rigidity, with particular emphasis on its relationship with conformal removability.
\end{abstract}

\maketitle

\section{Introduction}

Throughout this article, we will denote the complex plane by $\mC$, the Riemann sphere by $\hC:= \mC \cup \{\infty\}$ and the open unit disk by $\mD:= \{z \in \mC: |z|<1\}$. A \textit{domain} in $\hC$ is a connected open set. We say that two domains $U,V \subset \hC$ are \textit{conformally equivalent} if there exists a conformal map $f:U \to V$.

It is a remarkable fact, first formulated heuristically by Riemann in 1851 and proved rigorously later, that every simply connected proper domain in $\mC$ is conformally equivalent to the unit disk.

\begin{thm}[Riemann Mapping Theorem]
Every simply connected domain $\Omega \subsetneq \mC$ is conformally equivalent to $\mD$. Moreover, a conformal map $f:\Omega \to \mD$ is unique up to composition with a M\"{o}bius automorphism of $\mD$.
\end{thm}

\begin{figure}[ht]
\centering
\begin{tikzpicture}[x=.85cm,y=.85cm]
  \fill[domainpink] (0,0) rectangle (10,6.8);
  \begin{scope}
    \clip (0,0) rectangle (10,6.8);
    \foreach \k in {0,...,145} {
      \pgfmathsetmacro{\a}{mod(137.508*mod(\k,31),360)+floor(\k/31)*3}
      \pgfmathsetmacro{\r}{.76+.0088*mod(37*\k,69)}
      \fill ({6.55+\r*cos(\a)},{5.13+\r*sin(\a)}) circle[radius=.012];
    }
    \foreach \k in {0,...,28} {
      \pgfmathsetmacro{\a}{mod(137.508*mod(\k,31),360)}
      \pgfmathsetmacro{\r}{.08+.008*\k}
      \fill ({8.42+\r*cos(\a)},{1.72+\r*sin(\a)}) circle[radius=.013];
    }
    \foreach \p in {(8.17,2.45),(8.30,2.36),(8.41,2.53),
       (9.14,1.90),(9.25,2.02),(9.33,1.86),
       (9.53,3.00),(9.58,2.91),(8.96,2.68)}
       \fill \p circle[radius=.025];
  \end{scope}

  \draw[hole] (2.50,4.05) circle[radius=1.70];
  \draw[hole] (6.55,5.13) circle[radius=.60];
  \draw[hole] (8.92,3.70) circle[radius=.43];
  \draw[hole] (9.65,3.35) circle[radius=.17];
  \draw[hole] (8.60,3.14) circle[radius=.10];
  \foreach \x/\y/\r in {
       3.70/2.00/.16,4.12/2.13/.20,4.55/2.30/.12,
       4.55/1.55/.09,4.07/1.63/.11,
       5.05/1.30/.34,5.92/.92/.47,6.39/1.62/.26,
       5.43/1.86/.15,4.58/.83/.14,4.87/.48/.16,
       5.42/.37/.10,6.48/.39/.09,6.72/.69/.08,
       6.81/1.04/.09,6.77/1.37/.08}
       \draw[hole] (\x,\y) circle[radius=\r];
\end{tikzpicture}
 \caption{A circle domain. Its boundary can have isolated circles, circles converging to points, points converging to circles, Cantor sets, etc.}
 \label{fig:circle-domain}
\end{figure}
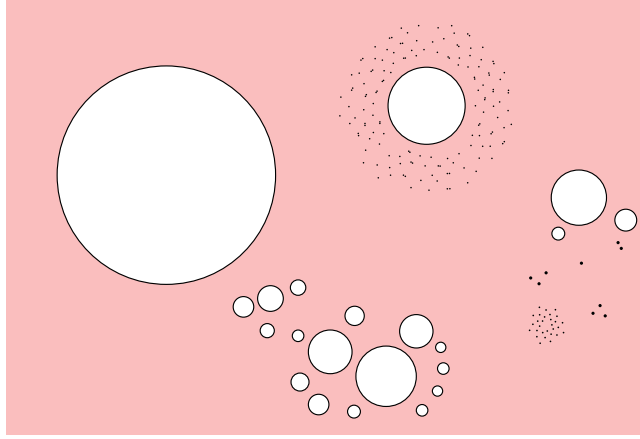

The unit disk $\mD$ is the simplest example of a \textit{circle domain}, that is, a domain $\Omega \subset \hC$ for which every connected component of $\partial \Omega$ is either a circle or a point. See Figure \ref{fig:circle-domain}.


The next result generalizes the Riemann mapping theorem to the setting of \textit{finitely connected domains}, that is, domains with finitely many complementary components.

\begin{thm}[Koebe, 1920 {\cite{Koe20}}]
\label{thm:koebe}
Every finitely connected domain $\Omega \subset \hC$ is conformally equivalent to a circle domain $D$. Moreover, a conformal map of $\Omega$ onto a circle domain is unique up to composition with a M\"{o}bius transformation.
\end{thm}

Koebe had earlier conjectured that conformal equivalence with a circle domain should hold for \textit{all} domains in $\hC$.

\begin{con}[Koebe's conjecture, 1908 {\cite{Koe08}}]
Every domain in $\hC$ is conformally equivalent to a circle domain.
\end{con}

This deep problem remains open after more than a hundred years, despite major advances by He--Schramm \cite{HS93} and many others. For more information on Koebe's conjecture and recent progress, we refer the reader to the survey article \cite{Nta26Survey}.

Koebe's conjecture concerns the \textit{existence} of a conformal map $f:\Omega \to D$ from a given domain $\Omega \subset \hC$ onto a circle domain $D$. Here we focus instead on the corresponding \textit{uniqueness} question. Composing $f$ with a M\"{o}bius transformation yields another conformal map $g:\Omega \to D'$, where $D'$ is also a circle domain. Must every conformal map from $\Omega$ onto a circle domain arise in this way? This leads to the notion of conformal rigidity: a circle domain $D \subset \hC$ is \textit{conformally rigid} if every conformal map from $D$ onto another circle domain is the restriction of a M\"{o}bius transformation.

\begin{figure}[ht]
\centering
\begin{tikzpicture}[x=.85cm,y=.85cm]
  \fill[domainpink] (0,0) rectangle (5.65,5.1);
  \draw[hole] (0.55,2.00) .. controls (.55,2.24) and (.90,2.50) .. (1.15,2.60)
    .. controls (1.42,2.73) and (1.39,3.02) .. (1.39,3.52)
    .. controls (1.39,3.81) and (1.58,3.93) .. (1.84,3.93)
    -- (2.32,3.93) .. controls (2.59,3.93) and (2.71,3.80) .. (2.70,3.55)
    -- (2.66,2.67) .. controls (2.65,2.49) and (2.51,2.38) .. (2.26,2.29)
    -- (.73,1.76) .. controls (.57,1.71) and (.47,1.80) .. (.55,2.00) -- cycle;
  \draw[hole] (3.48,3.63) .. controls (3.43,3.94) and (3.60,4.03) .. (3.77,3.79)
    -- (4.03,3.46) .. controls (4.22,3.24) and (4.40,3.37) .. (4.70,3.53)
    .. controls (4.99,3.70) and (5.08,3.60) .. (5.04,3.29)
    -- (4.99,2.74) .. controls (4.98,2.57) and (4.79,2.46) .. (4.56,2.44)
    -- (3.59,2.41) .. controls (3.37,2.41) and (3.33,2.57) .. (3.38,2.82) -- cycle;
  \draw[hole] (2.31,.62) .. controls (2.14,.98) and (2.10,1.15) .. (2.34,1.40)
    -- (2.79,1.85) .. controls (2.94,2.00) and (3.08,1.91) .. (3.27,1.77)
    -- (3.73,1.43) .. controls (3.86,1.34) and (3.99,1.37) .. (4.15,1.49)
    -- (4.37,1.68) -- (4.35,.64) .. controls (4.35,.47) and (4.24,.42) .. (4.05,.43)
    -- (2.66,.45) .. controls (2.47,.45) and (2.39,.48) .. (2.31,.62) -- cycle;
  \draw[hole] (1.34,1.76) .. controls (1.40,1.98) and (1.61,2.00) .. (1.78,1.92)
    -- (1.94,1.80) .. controls (2.02,1.72) and (1.93,1.58) .. (1.84,1.49)
    -- (1.66,1.37) .. controls (1.49,1.28) and (1.23,1.48) .. (1.34,1.76) -- cycle;
  \draw[hole] (3.05,2.02) -- (4.86,2.02) -- (4.77,2.14)
     .. controls (4.72,2.21) and (4.65,2.22) .. (4.56,2.22)
     -- (3.25,2.22) -- cycle;
  \draw[hole] (3.02,2.29) -- (4.84,2.29) -- (4.70,2.36) -- (3.17,2.36) -- cycle;
  \draw[hole] (4.66,.45) -- (4.85,.44) -- (4.80,1.72)
      .. controls (4.77,1.82) and (4.68,1.85) .. (4.65,1.73) -- cycle;
  \draw[hole] (4.96,.45) -- (5.07,.45) -- (5.04,1.76) -- (4.96,1.78) -- cycle;
  \foreach \x/\y/\r in {2.96/2.67/.10,1.64/1.03/.07}
    \draw[hole] (\x,\y) circle[radius=\r];
  \foreach \x/\y/\dx/\dy in {
     .42/.42/-.13/.72,.54/.51/-.12/.58,.75/.39/-.18/.75,
     .91/.56/-.12/.62,1.08/.52/-.18/.73,1.37/.41/-.12/.62,
     1.72/.27/-.13/.54,1.91/.23/-.11/.61}
      \draw[line width=.45pt] (\x,\y) -- ++(\dx,\dy);

  \fill[domainpink] (7.60,3.05) rectangle (11.75,5.10);
  \fill[domainpink] (7.60,0) rectangle (11.75,2.05);
  \foreach \x/\y/\r in {
    9.65/3.92/.59,8.65/4.47/.40,10.70/4.55/.31,
    10.76/3.54/.24,8.55/3.50/.20,9.11/3.24/.18,
    8.05/4.16/.10,8.07/3.80/.11,9.02/4.94/.13,
    11.33/3.97/.12,11.35/3.65/.09,10.14/3.23/.10,
    10.44/3.23/.09,9.57/3.19/.08,
    10.52/.75/.60,9.10/1.31/.40,11.29/1.47/.31,
    8.15/1.68/.11,8.17/1.34/.10,8.79/.68/.21,
    9.60/.34/.16,11.17/.20/.08,11.47/.85/.11,
    11.46/.56/.08,9.33/.81/.07,9.53/.85/.05}
    \draw[hole] (\x,\y) circle[radius=\r];
  \foreach \x/\y in {
    8.34/3.17,8.45/3.14,8.56/3.13,8.67/3.15,
    8.27/1.10,8.31/1.02,9.72/.91,9.72/.68}
    \fill (\x,\y) circle[radius=.018];
  \draw[transformation] (5.98,2.72) -- (7.23,3.93);
  \draw[transformation] (5.98,2.35) -- (7.23,1.12);
  \draw[transformation] (9.68,2.88) -- (9.68,2.23);
\end{tikzpicture}
\caption{If these two conformal maps from the domain on the left to the circle domains on the right do not differ by a M\"{o}bius transformation, then neither target circle domain is conformally rigid.}
 \label{fig:uniqueness}
\end{figure}
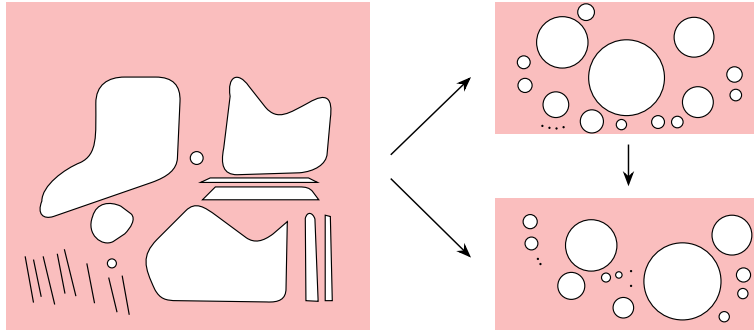

For instance, finitely connected circle domains are conformally rigid by Theorem \ref{thm:koebe}. This was generalized by He and Schramm in \cite{HS93}, who proved that every circle domain with at most countably many boundary components is conformally rigid. They later showed that circle domains whose boundaries have $\sigma$-finite length are conformally rigid \cite{HS94}.

On the other hand, it is well known that not every circle domain is conformally rigid. For instance, if $K \subset \mC$ is a Cantor set of positive area, one can use the theory of quasiconformal mappings to construct a homeomorphism $f:\hC \to \hC$ that is conformal on $\hC \setminus K$ but is not a M\"{o}bius transformation. Then $f$ maps the circle domain $\Omega:=\hC \setminus K$ onto another circle domain, but is not the restriction of a M\"{o}bius transformation, showing that $\Omega$ is not conformally rigid. See e.g. \cite[Proposition 4.4]{You15} and the remark after.

The existence of non-rigid circle domains is considered a major obstacle to resolving Koebe's conjecture. A better understanding of conformal rigidity may therefore shed light on the conjecture.

\section{Conformal rigidity and removability}

Recall from the Introduction that if $K \subset \mC$ is a Cantor set of positive area, then $\Omega:=\hC \setminus K$ is a circle domain that is not conformally rigid. This follows from the fact that $K$ is not conformally removable. A compact set $E \subset \hC$ is said to be \textit{conformally removable} if every homeomorphism of $\hC$ that is conformal on $\hC \setminus E$ is a M\"{o}bius transformation. Examples of conformally removable sets include sets of $\sigma$-finite length, by a theorem of Besicovitch \cite{Bes31}, and boundaries of uniform domains, by a theorem of Jones \cite{Jon95}. Examples of non-removable sets include sets of positive area, the Sierpinski gasket \cite{Nta19} and Sierpinski carpets \cite{Nta21}. We refer the reader to \cite{You15} and \cite{Bis26} for further information and open problems concerning conformal removability.

Regarding the relationship between conformal rigidity and conformal removability, He and Schramm posed the following conjecture.

\begin{con}[Rigidity conjecture, 1994 {\cite{HS94}}]
\label{con:rigidity}
A circle domain $\Omega$ is conformally rigid if and only if its boundary $\partial \Omega$ is conformally removable.
\end{con}

Conjecture~\ref{con:rigidity} is known to hold for several classes of circle domains, including those satisfying a certain condition on the quasihyperbolic distance. This was proved by the authors in \cite{NY20} and subsequently generalized by the first author in \cite{Nta23c}. More recently, Rajala disproved the conjecture by constructing a conformally rigid circle domain $\Omega$ whose boundary is not conformally removable \cite{Raj25}. The table below summarizes these developments, which we discuss in later sections.

\begin{table}[htbp]
\centering
\small
\begin{tabular}{@{}>{\raggedright\arraybackslash}p{0.36\textwidth}>{\raggedright\arraybackslash}p{0.30\textwidth}>{\raggedright\arraybackslash}p{0.27\textwidth}@{}}
\hline
Types of circle domains & $\Omega$ rigid & $\partial\Omega$ removable \\
\hline
Simply connected & \greencheck & \greencheck\,Morera \\
Finitely connected & \greencheck\,Koebe \cite{Koe20} & \greencheck\,Morera \\
Countably connected & \greencheck\,He--Schramm \cite{HS93} & \greencheck\,Besicovitch \cite{Bes31} \\
Boundary has $\sigma$-finite length & \greencheck\,He--Schramm \cite{HS94} & \greencheck\,Besicovitch \cite{Bes31} \\
Uniform & \greencheck\,Ntalampekos--Younsi \cite{NY20} & \greencheck\,Jones \cite{Jon95} \\
Quasihyperbolic distance in $L^2$ & \greencheck\,Ntalampekos--Younsi \cite{NY20} & \greencheck\,Jones--Smirnov \cite{JS00} \\
Boundary is CNED & \greencheck\,Ntalampekos \cite{Nta23c} & \greencheck\,Ntalampekos \cite{Nmetric} \\
Boundary has positive area & \redcross\,Sibner \cite{Sib68} & \redcross\, \\
Rajala's example & \greencheck\,Rajala \cite{Raj25} & \redcross\,Wu \cite{Wu98,Raj25} \\
\hline
\end{tabular}
\caption{Results related to the conformal rigidity of circle domains and the removability of their boundaries.}
\label{tab:rigidity-removability}
\end{table}

\section{A quasihyperbolic condition for rigidity}

For a domain $D\subsetneq \mC$ and $z\in D$ we set $\delta_D(z)=\dist(z,\partial D)$. For $x,y\in D$ we define their \textit{quasihyperbolic distance} in $D$ by
\begin{equation}
k_D(x,y)=\inf_\gamma\int_\gamma\frac{|dz|}{\delta_D(z)},
\label{eq:qhmetric}
\end{equation}
where the infimum ranges over all rectifiable paths in $D$ from $x$ to $y$. Paths that remain close to the boundary are expensive in this metric. Thus, $k_D$ detects narrow passages and poorly accessible portions of the domain.

The quasihyperbolic distance can be used to obtain a sufficient condition for the conformal rigidity of circle domains, as the following theorem shows.

\begin{thm}[{\cite{NY20}}]
\label{thm:qh}
Let $\Omega\subset\hC$ be a circle domain with $\infty\in\Omega$. Let $R>0$ such that $\hC\setminus \Omega\subset B(0,R)$ and set $D=\Omega\cap B(0,R)$. If, for some $x_0\in D$,
\begin{equation}
\int_D k_D(x,x_0)^2\,dA(x)<\infty,
\label{eq:qhcondition}
\end{equation}
then $\Omega$ is conformally rigid.
\end{thm}
Here $dA$ denotes integration with respect to Lebesgue measure. The truncation of the original domain $\Omega$ avoids integration near infinity. Since $D$ has finite area, Condition \eqref{eq:qhcondition} is independent of the choice of the base point $x_0$.

Jones and Smirnov \cite{JS00} proved removability under this quasihyperbolic integrability condition. Thus Theorem~\ref{thm:qh} establishes the rigidity conjecture (Conjecture~\ref{con:rigidity}) for this class. The condition includes uniform, John, and H\"older domains, giving the following consequence.
\begin{cor}
Uniform, John, and H\"older circle domains are conformally rigid.
\end{cor}
A bounded John domain $D$ has a base point $x_0$ such that every $x\in D$ can be joined to $x_0$ by an arc-length parametrized curve $\gamma\colon[0,L]\to D$ satisfying $\gamma(0)=x, \gamma(L)=x_0$ and $\delta_D(\gamma(t))\ge ct$, with constant $c>0$ independent of $x$. This condition rules out arbitrarily long narrow corridors in the domain $D$. One standard quasihyperbolic formulation of the H\"older domain condition is
\begin{equation}
k_D(x,x_0)\le A\log\frac{\delta_D(x_0)}{\delta_D(x)}+B,
\label{eq:holder}
\end{equation}
for fixed constants $A,B>0$ and for a base point $x_0\in D$. Every John domain is also a H\"older domain. The relevant integrability properties of the quasihyperbolic distance in these domains are discussed in \cite{SS}. 

These hypotheses constrain how one approaches the boundary from the domain. They neither require a rectifiable boundary nor a bound on its total length. This explains how Theorem \ref{thm:qh} can extend beyond the $\sigma$-finite-length theorem of He and Schramm from \cite{HS94}. Conversely, a measure-theoretic hypothesis on the boundary does not by itself rule out severe bottlenecks inside the domain. Hence, our result does not imply the He--Schramm theorem.

We discuss the idea of the proof of Theorem \ref{thm:qh}. If the boundary of a domain has $\sigma$-finite length, as in the He--Schramm theorem, then almost every line or circle meets it in at most countably many points. This crucial property gives absolute continuity estimates, which are  used in combination with modulus arguments in the proof of He and Schramm. Under the quasihyperbolic condition, the same intersection statement can fail. The replacement is to deform paths slightly, introducing \emph{detours} of paths that have a complicated intersection with the boundary. Namely, if a line segment intersects the boundary of the domain at uncountably many points, then we deform it slightly by using quasihyperbolic geodesics to a path that intersects the boundary at finitely many points. Condition \eqref{eq:qhcondition} ensures the existence of detours that have small cost. We now briefly describe the main steps of the proof of Theorem \ref{thm:qh}.

\emph{Boundary extension.} Suppose that $\Omega$ is a circle domain satisfying \eqref{eq:qhcondition} and let $f$ be a conformal map from $\Omega$ onto another circle domain $\Omega^*$. One first uses modulus estimates and detours to prove that point boundary components of $\Omega$ correspond exactly to point boundary components of $\Omega^*$ under $f$. Additional estimates show that $f$ extends continuously and injectively to $\overline \Omega$. This produces a homeomorphism $\overline\Omega\to\overline{\Omega^*}$.

\emph{Reflection and analytic control.} One repeatedly reflects across matching circles to obtain a homeomorphism $F:\hC \to \hC$. This extension conjugates the group generated by reflections in the boundary circles of $\Omega$, i.e., the \textit{Schottky group} of $\Omega$, to the Schottky group of $\Omega^*$. The detour estimates are used again to obtain modulus bounds for $F$, and hence quasiconformality. Crucially, the quasiconformal dilatation of $F$ is bounded by some $K_0$ that is independent of the choice of the conformal map $f$ and its target circle domain $\Omega^*$.

\emph{Quasiconformal deformation.} If $F$ is not a M\"obius transformation, then it has a non-trivial Beltrami coefficient, which can be modified in a way compatible with the Schottky group. Solving the corresponding Beltrami equation yields a quasiconformal homeomorphism of the sphere that maps $\Omega$ conformally onto another circle domain and whose quasiconformal dilatation exceeds $K_0$. This is a contradiction, hence the original conformal map $f:\Omega \to \Omega^*$ must be a M\"obius transformation.

\section{A general condition for rigidity}

In an effort to find a necessary and sufficient condition for the removability of a compact set, the first author introduced the notion of a CNED set, that is, \textit{countably negligible for extremal distances}. In order to give the definition, we introduce the notion of conformal modulus of a curve family in the Riemann sphere $\hC$. We equip $\hC$ with the spherical metric $\sigma$ and spherical measure $\Sigma$. Let $\Gamma$ be a family of curves in $\hC$. A Borel function $\rho\colon \hC\to [0,\infty]$ is \textit{admissible} for $\Gamma$ if
$$\int_\gamma\rho\, ds\geq1$$
for every rectifiable curve $\gamma\in \Gamma$. The \textit{conformal modulus} of $\Gamma$ is defined as
$$\Mod\Gamma = \inf_{\rho}\int_{\hC} \rho^2\, d\Sigma,$$
where the infimum is taken over all admissible functions $\rho$ for $\Gamma$. If there are no such admissible functions, then we set $\Mod\Gamma=\infty$.

Let $E\subset \hC$ be a set and $\mathcal F_{\sigma}(E)$ be the family of curves in $\hC$ that intersect $E$ at countably many points, not counting multiplicity. For two sets $F_1,F_2\subset \hC$ we denote by $\Gamma(F_1,F_2)$ the family of curves $\gamma\colon[a,b]\to \hC$ with $\gamma(a)\in F_1$ and $\gamma(b)\in F_2$. We say that a set $E\subset \hC$ is CNED if for every pair of non-empty, disjoint continua $F_1,F_2\subset \hC$ we have
$$\Mod\Gamma(F_1,F_2)=\Mod(\Gamma(F_1,F_2)\cap \mathcal F_\sigma(E)).$$
In other words, the modulus of the curve family $\Gamma(F_1,F_2)$ is not affected when we restrict to the subfamily of curves intersecting $E$ at countably many points.

This class of sets was introduced in \cite{Nmetric} and it was shown that compact CNED sets are conformally removable. Moreover, the CNED condition is conjectured to be necessary and sufficient for the conformal removability of a compact set. The properties of these sets were studied in depth in \cite{NCNED} and it was shown that this class includes sets of $\sigma$-finite length and boundaries of domains satisfying the quasihyperbolic condition of Theorem \ref{thm:qh}. A natural problem was to test this class of sets against the rigidity conjecture. The first author proved in \cite{Nta23c} that circle domains with CNED boundary are indeed conformally rigid.

\begin{thm}[\cite{Nta23c}]\label{thm:cned_rigid}
    A circle domain is conformally rigid if every compact subset of its point boundary components is CNED.
\end{thm}

This result unifies several different approaches to the problem of conformal rigidity. Moreover, it provides evidence that the conformal removability of the boundary of a circle domain implies the rigidity of the circle domain, as expected by Conjecture \ref{con:rigidity}.

We now discuss some ideas from the proof of Theorem \ref{thm:cned_rigid}, which is more conceptual than the proofs of earlier rigidity results. One of the main tools is a new metric characterization of (quasi)conformality due to the first author, which we describe below. For an open set $A\subset \hC$, define its \textit{eccentricity} by
\[
E(A)=\inf\{M\ge1:B\subset A\subset MB
\text{ for some ball }B\}.
\]
The ball $B$ is not required to be centered at a prescribed point. For a homeomorphism $F$, its \emph{eccentric distortion} $E_F(x)$ is the infimum of $H\ge1$ for which there exist open neighborhoods $A_n$ of $x$, shrinking to $x$ as $n\to\infty$, such that
\[
E(A_n)\le H \quad \text{and}\quad E(F(A_n))\le H
\]
for each $n\in \mathbb N$.

\begin{thm}[Metric criterion {\cite{Nmetric}}]
\label{thm:metric}
Let $F\colon U\to V$ be an orientation-preserving homeomorphism between open sets $U,V\subset \hC$, and let $E\subset \hC$ be a CNED set. If there exists $H\geq 1$ such that $E_F(x)\le H$ for all $x\in U\setminus E$, then $F$ is quasiconformal. If, in addition, $E_F(x)=1$ for almost every $x\in U$, then $F$ is conformal.
\end{thm}

We now discuss the main steps of the proof of Theorem \ref{thm:cned_rigid}, which is based on Theorem \ref{thm:metric}.

\textit{Extension and reflection.} Suppose that $f\colon \Omega\to \Omega^*$ is a conformal map between circle domains in $\hC$ such that $\Omega$ satisfies the condition of Theorem \ref{thm:cned_rigid}. The first step, as in the proof of Theorem \ref{thm:qh}, is to show that $f$ extends to a homeomorphism between the closures of the domains. This is achieved by using the CNED condition to modify curves whose intersection with $\partial \Omega$ is uncountable to curves with countable intersections, in the same way detours are used in the proof of Theorem \ref{thm:qh}. We then construct its reflection-compatible extension $F$ that is a homeomorphism of $\hC$, as in Theorem \ref{thm:qh}. Our goal is to show that $F$ is conformal. 

\textit{Conformality at interior points.} On each reflected copy of the original domain $\Omega$, the map $F$ is conformal and $E_F=1$, by the Schwarz reflection construction. At a residual point enclosed by nested reflected disks, those disks and their images shrink to the corresponding points, again giving $E_F=1$.

\textit{Negligibility of a large part of the boundary.}
The reflected boundary points require a more careful decomposition. Choose thin open neighborhoods of the complementary disks of $\Omega$ whose images under $F$ are also thin neighborhoods of the corresponding target disks in the complement of $\Omega^*$. The portions of the boundary of $\Omega$ outside these neighborhoods form compact subsets $K_n$ of the point components of $\partial \Omega$. By the main assumption of Theorem \ref{thm:cned_rigid}, each $K_n$ is a CNED set. Their reflected copies, together with the reflected circles, form a CNED set. This is because countable unions of compact CNED sets are CNED, as shown in \cite{NCNED}.

\textit{Eccentricity at point components near circles.}
The remaining point components of $\partial \Omega$ lie arbitrarily close to disks in the chosen neighborhoods. The neighborhoods and their images can be chosen with eccentricities tending to one, yielding $E_F=1$ at these remaining points as well.  

\textit{Metric characterization of conformality.}
By the above, we have $E_F=1$ outside a CNED set. The metric criterion of Theorem \ref{thm:metric} first supplies quasiconformal regularity of $F$. The exceptional set has area zero, so the last part of the theorem makes $F$ conformal, hence M\"obius. 

\section{Rigid domains with non-removable boundaries}

The preceding criteria support a close relationship between rigidity and removability, as predicted by Conjecture \ref{con:rigidity}. However, Rajala recently made a quite surprising discovery. 

\begin{thm}[Rajala {\cite{Raj25}}]
\label{thm:rajala}
There exists a conformally rigid circle domain whose boundary is not conformally removable.
\end{thm}
The boundary of this circle domain contains a non-removable Cantor set. Rajala's example disproves one direction of the rigidity conjecture. A key ingredient is the following theorem of Wu \cite{Wu98}.
\begin{thm}[Wu \cite{Wu98}]
Let $E,F\subset\mR$ be Cantor sets and $[a,b]\subset \mR$. If
\begin{equation}
\caplog([a,b]\setminus F)<\caplog([a,b]),
\label{eq:wu}
\end{equation}
where $\caplog$ is logarithmic capacity, then $E\times F$ is conformally non-removable.
\end{thm}
Condition~\eqref{eq:wu} is a potential-theoretic thickness condition on $F$. One can retain it while choosing the horizontal Cantor set $E$ very thin. In Rajala's construction $E$ has length zero, so $K=E\times F$ has area zero by Fubini's theorem. Thus the non-removability of $E\times F$, provided by Wu's theorem, is non-trivial.

\begin{figure}
    \centering
    \includegraphics[width=0.8\linewidth]{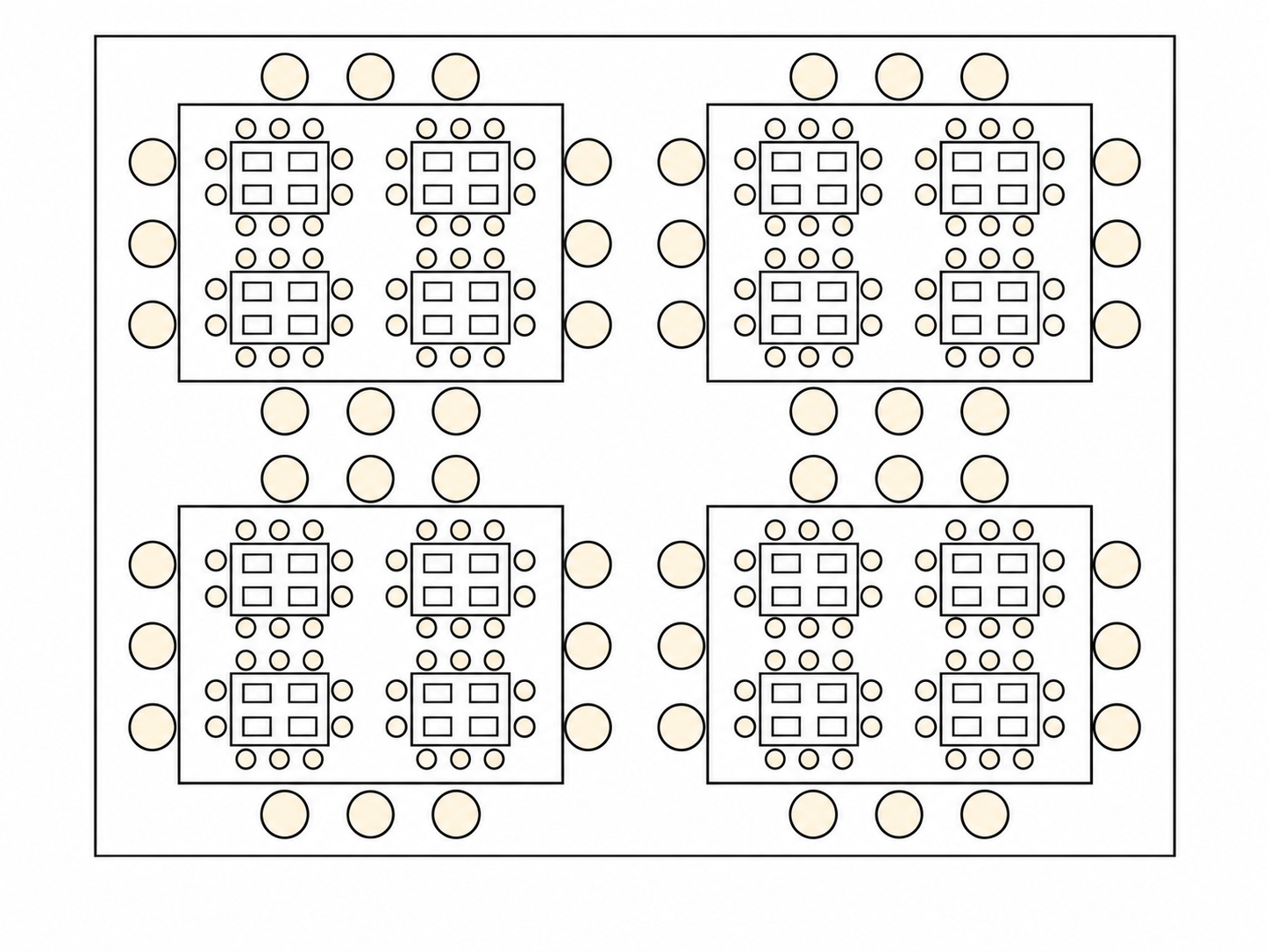}
    \caption{Illustration of Rajala's example. Countably many isolated disks $D_j$ accumulate at the points of the product Cantor set $K=E\times F$.}
    \label{fig:rajala}
\end{figure}

Rajala chooses pairwise disjoint closed disks $D_j$ that are disjoint from $K$ and accumulate at $K$, and then considers the domain
\begin{equation}
\Omega=\hC\setminus\bigg(K\cup\bigcup_j D_j\bigg).
\label{eq:rajala}
\end{equation}
The disks $D_j$ are placed at successively smaller scales alongside rectangles defining the product Cantor set. The construction balances two requirements: the vertical Cantor set must remain thick enough for Wu's theorem, while the horizontal scales and neighboring disks must enforce strong modulus estimates. See Figure \ref{fig:rajala} for an illustration of the construction.

The non-removability of $\partial \Omega$ immediately follows from the non-removability of $K\subset\partial\Omega$. Indeed, a non-M\"obius sphere homeomorphism conformal outside $K$ is certainly conformal outside the larger set $\partial\Omega$. The rigidity of $\Omega$ requires a separate argument. The neighboring disks force a conformal map $f$ from $\Omega$ onto an arbitrary circle domain to obey quantitative distortion bounds near $K$. Modulus estimates establish boundary extension of $f$. As usual, we then consider a reflection-compatible extension of $f$ to a homeomorphism $G$ of $\hC$. Furthermore, modulus estimates control the relative positions and sizes of the images of the disks $D_j$ under $G$. These bounds give uniformly bounded eccentric distortion of $G$ at points of $K$ as well as its reflected copies. Away from them the extension has eccentric distortion one, including at residual points surrounded by nested round disks. Compare to the proof outline of Theorem \ref{thm:cned_rigid}.

The metric criterion of Theorem \ref{thm:metric}, now with no exceptional CNED set needed, gives global quasiconformality of $G$. The reflected copies of $K$ have area zero, so the eccentric distortion of $G$ equals one almost everywhere. The extension is therefore conformal. This proves rigidity even though $K$ itself cannot be CNED, since compact CNED sets are removable.

The example shows that the circles of the boundary of a circle domain have an essential role. The complement $\hC\setminus K$ is non-rigid, while removing the carefully arranged disks produces a rigid circle domain. Thus rigidity is not determined solely by the point components of $\partial \Omega$ and is not monotone under inclusion of domain boundaries. A homeomorphism $g\colon \hC\to \hC$ witnessing the non-removability of $K$ (so $g$ is conformal in $\hC\setminus K$ but not M\"obius) cannot map all the disks $D_j$ to round disks, since that would contradict the rigidity of $\Omega$. However, the requirement that the images of $D_j$ are round disks (as in the conformal rigidity problem) imposes constraints that are absent from the removability problem and that force $g$ to be a M\"obius transformation.

\section{Remaining implication in the rigidity conjecture}
\begin{figure}
    \centering
    \includegraphics[width=0.6\linewidth]{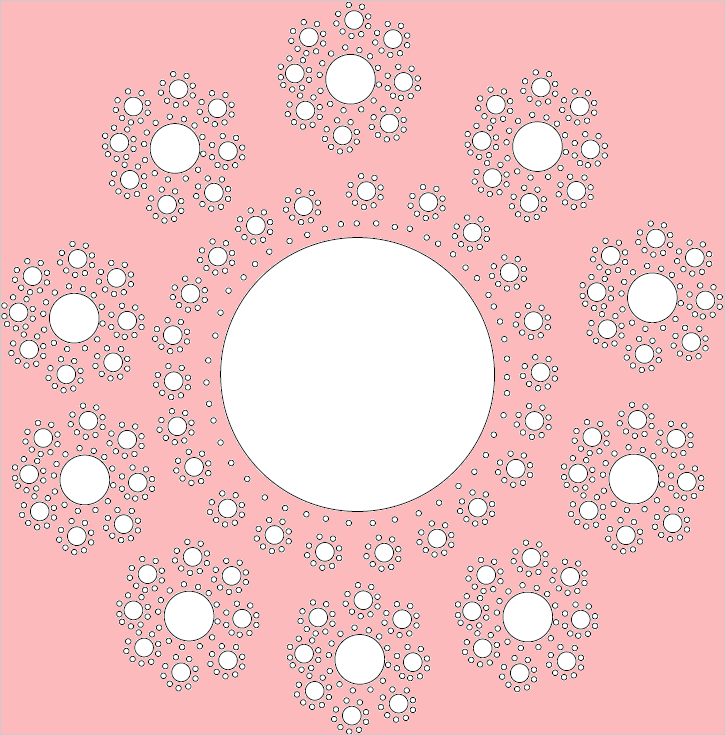}
    \caption{A circle domain in which every point of each boundary circle is an accumulation point of infinitely many distinct boundary circles.}
    \label{fig:sierpinski_type}
\end{figure}

\begin{ques}\label{q:remaining}
If a circle domain has conformally removable boundary, must it be conformally rigid?
\end{ques}
This remains open even when the entire boundary is a Cantor set. In that case one implication is immediate from the definitions: rigidity of $\hC\setminus K$ implies removability of $K$. For the converse there is a useful reduction. A conformal homeomorphism between the complements of two compact totally disconnected sets extends to a homeomorphism of the sphere. Consequently, if $K$ is removable and
\[
f\colon\hC\setminus K\longrightarrow\Omega^*
\]
is a non-M\"obius conformal map onto a circle domain $\Omega^*$, then $\partial\Omega^*$ must contain a circle. That is, a point component of the source must open into a non-degenerate component of the target. It is this possibility that prevents removability from settling the question. The authors proved a stronger restriction \cite[Theorem~9.1]{NY20}: in such a hypothetical counterexample, every point on every boundary circle of $\Omega^*$ must be an accumulation point of infinitely many distinct boundary circles. Thus an isolated target circle cannot occur. The target must exhibit a Sierpi\'nski-type accumulation pattern, as in Figure \ref{fig:sierpinski_type}. This narrows the geometry of a possible counterexample, but does not rule it out.

A further conjecture places Question~\ref{q:remaining} within the general theory of removability.
\begin{con}[CNED characterization {\cite{NCNED}}]
A compact planar set is conformally removable if and only if it is CNED.
\label{con:cned}
\end{con}
Only the implication from removability to CNED is open. If it holds, a removable circle domain boundary is CNED, and Theorem~\ref{thm:cned_rigid} gives rigidity.  Rajala's example is compatible with this conjecture: its non-removable Cantor subset fails the CNED condition, while surrounding circles enforce rigidity by a different mechanism. Summarizing, a proof of Conjecture \ref{con:cned} would provide a positive answer to Question \ref{q:remaining} as well as a proof of Conjecture \ref{con:rigidity} for circle domains with totally disconnected boundaries.

We conclude with a general criterion for the rigidity of circle domains, which may also provide an approach to Question~\ref{q:remaining}.

\begin{thm}[\cite{NtalampekosRajala:exhaustion}]
Let $\Omega\subset \hC$ be a circle domain and $g$ be a conformal map from $\Omega$ onto another circle domain. Suppose that there exists a curve family $\Gamma_0$ in $\hC$ with $\Mod\Gamma_0=0$ such that for every $\gamma\notin \Gamma_0$ with endpoints $z,w\in \Omega$ we have
$$\sigma(g(z),g(w))\leq \ell(g\circ \gamma|_{\gamma^{-1}(\Omega)})+ \sum_{q\cap |\gamma|\neq \emptyset} \diam \hat g(q),$$
where the sum ranges over all complementary components $q$ of $\Omega$ that intersect the curve $\gamma$ and $\hat g(q)$ denotes the complementary component of $g(\Omega)$ corresponding to $q$. Then $g$ is the restriction of a M\"obius transformation.
\end{thm}

\end{document}